\documentclass[12pt]{amsart}

\usepackage{amssymb}
\usepackage{graphicx}
\usepackage{epsf}
\usepackage{amsmath}
\usepackage[latin1]{inputenc}
\usepackage{amsfonts}
\usepackage{color}

\numberwithin{equation}{section}
\newtheorem{theo}{Theorem}[section]
\newtheorem{prop}[theo]{Proposition}

\newtheorem{rema}[theo]{Remark}

\newcommand{\pref}[1]{(\ref{#1})}

\def\ie{{\it i.e. }}

\def\ti{$\sim$}

\def\C{{\mathbb C}}
\def\N{{\mathbb N}}

\def\Q{{\mathbb Q}}

\def\I{{\mathcal I}}

\def\D{{\mathcal D}}

\def\T{{\mathcal T}}

\def\CC{{\mathcal C}}

\def\QQ{{\bf Q}}

\def\shuffle{{\,\raise
1pt\hbox{$\scriptscriptstyle\cup{\mskip-4mu}\cup$}\,}}

\title{Multivariate Fuss-Catalan numbers}

\author{J.-C.~Aval}\address[Jean-Christophe Aval]{LaBRI\\ Universit\'e Bordeaux 1\\ 351 cours
de 
 la Lib\'eration\\ 33405 Talence cedex\\ FRANCE}
\email{aval@labri.fr}
\urladdr{http://www.labri.fr/\ti aval/}

\date{\today}

\begin{document} 

\begin{abstract} 
Catalan numbers $C(n)=\frac{1}{n+1}{2n\choose n}$ enumerate binary trees and Dyck paths. The distribution of paths with respect to their number $k$ of factors is given by ballot numbers $B(n,k)=\frac{n-k}{n+k}{n+k\choose n}$. These integers are known to satisfy simple recurrence, which may be visualised in a ``Catalan triangle'', a lower-triangular two-dimensional array. It is surprising that the extension of this construction to 3 dimensions generates integers $B_3(n,k,l)$ that give a 2-parameter distribution of $C_3(n)=\frac 1 {2n+1} {3n\choose n}$, which may be called order-3 Fuss-Catalan numbers, and enumerate ternary trees. The aim of this paper is a study of these integers $B_3(n,k,l)$. We obtain an explicit formula and a description in terms of trees and paths. Finally, we extend our construction to $p$-dimensional arrays, and in this case we obtain a $(p-1)$-parameter distribution of $C_p(n)=\frac 1 {(p-1)n+1} {pn\choose n}$, the number of $p$-ary trees.
\end{abstract}   

\maketitle

\section{Catalan triangle, binary trees, and Dyck paths}

We recall in this section well-known results about Catalan numbers and ballot numbers.

The {\em Catalan numbers} 
$$C(n)=\frac{1}{n+1}{2n\choose n}$$
are integers that appear in many combinatorial problems. These numbers first appeared in Euler's work as the number of triangulations of a polygon by mean of non-intersecting diagonals. Stanley \cite{stanley,stanweb} maintains a dynamic list of exercises related to Catalan numbers, including (at this date) 127 combinatorial interpretations. 

Closely related to Catalan numbers are {\em ballot numbers}. Their name is due to the fact that they are the solution of the so-called {\em ballot problem}: we consider an election between two candidates A and B, which respectively receive $a>b$ votes. The question is: what is the probability that during the counting of votes, A stays ahead of B? The answer will be given below, and we refer to \cite{bertrand} for Bertrand's first solution, and to \cite{andre} for Andr\'e's beautiful solution using the ``reflection principle''.

Since our goal here is different, we shall neither define ballot numbers by the previous statement, nor by their explicit formula, but we introduce integers $B(n,k)$ defined for a positive integer $n$ and a nonnegative integer $k$ by the following conditions:
\begin{itemize}
\item $B(1,0)=1$;
\item $\forall n> 1$ and $0\le k< n$, $B(n,k)=\sum_{i=0}^k B(n-1,i)$;
\item $\forall k\ge n$, $B(n,k)=0$.
\end{itemize}

Observe that the recursive formula in the second condition is equivalent to: 
\begin{equation}\label{recc}
B(n,k)=B(n-1,k)+B(n,k-1). 
\end{equation}

We shall present the $B(n,k)$'s by the following triangular representation (zero entries are omitted) where moving down increases $n$ and moving right increases $k$.

\vskip 0.2cm
$
\begin{array}{rrrrrr}
1&&&&&\cr
1&1&&&&\cr
1&2&2&&&\cr
1&3&5&5&&\cr
1&4&9&14&14&\cr
1&5&14&28&42&42\cr
\end{array}
$
\vskip 0.2cm

The crucial observation is that computing the horizontal sums of these integers give : $1,\ 2,\ 5,\ 14,\ 42,\ 132$. We recognize the first terms of the Catalan series, and this intuition will be settled in Proposition \ref{prop1}, after introducing combinatorial objects. 

\vskip 0.2cm
A {\em binary tree} is a tree in which every internal node has exactly 2 sons. The number of binary trees with $n$ internal nodes is given by the $n$-th Catalan number. The nodes of the following tree are labelled to explain a bijection described in the next paragraph: internal nodes are labelled by letters and external nodes by numbers.

\vskip 0.2cm
{\centerline{\epsffile{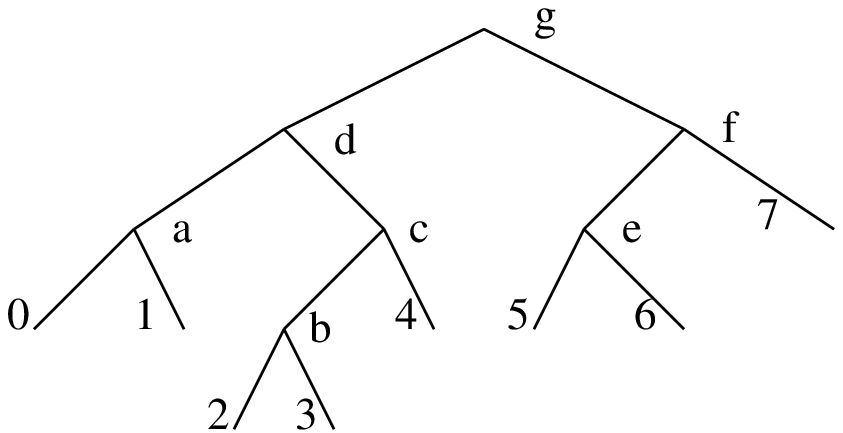}}}

A {\em Dyck path} is a path consisting of steps $(1,1)$ and $(1,-1)$, starting from $(0,0)$, ending at $(2n,0)$, and remaining in the half-plane $y\ge0$ (we shall sometimes say ``remaining above the horizontal axis'' in the same sense). The number of Dyck paths of length $2n$ is also given by the $n$-th catalan number. More precisely, the depth-first search of the tree gives a bijection between binary trees and Dyck paths: we associate to each external node (except the left-most one) a $(1,1)$ step and to each internal node a $(1,-1)$ step by searching recursively the left son, then the right son, then the root. As an example, we show below the Dyck path corresponding to the binary tree given above. The labels on steps correspond to those on the nodes of the tree.

{\centerline{\epsffile{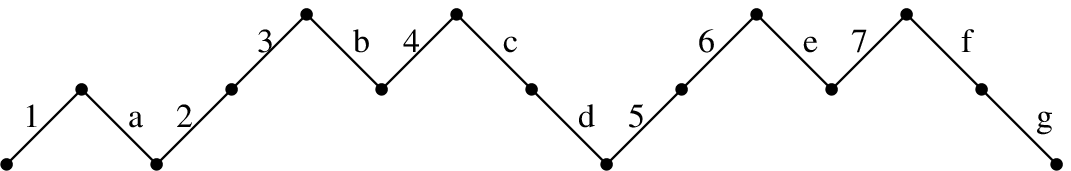}}}

An important parameter in our study will be the length of the right-most sequence of $(1,-1$) of the path. This parameter equals 2 in our example.
Observe that under the correspondence between paths and trees, this parameter corresponds to the length of the right-most string of right sons in the tree. We shall use the expressions {\em last down sequence} and {\it last right string}, for these parts of the path and of the tree.

Now we come to the announced result. It is well-known and simple, but is the starting point of our work.

\begin{prop}\label{prop1}
For $n$ a positive integer, we have the following equality:
$$\sum_{k=0}^{n-1} B(n,k)=C(n)=\frac 1 {n+1} {2n \choose n}.$$
\end{prop}

\begin{proof}
Let us denote by $\CC_{n,k}$ the set of Dyck paths of length $2n$ with a last down sequence of length equal to $n-k$.

We shall prove that $B(n,k)$ is the cardinality of $\CC_{n,k}$.

The proof is done recursively on $n$. If $n=1$, this is trivial. If $n>1$, let us suppose that $B(n-1,k)$ is the cardinality of $\CC_{n-1,k}$ for $0\le k<n-1$. Let us consider an element of $\CC_{n,k}$. If we erase the last step $(1,1)$ and the following step $(1,-1)$, we obtain a Dyck path of length $2(n-1)$, with a last decreasing sequence of length $n-1-l$ with $l\le k$. If we keep track of the integer $k$, we obtain a bijection between $\CC_{n,k}$ and $\cup_{l\le k}\CC_{n-1,l}$. We mention that this process is very similar to the ECO method \cite{eco}.

This is a combinatorial proof of Proposition \ref{prop1}.
\end{proof}

\noindent {\bf Remark 1.2.}
The integers $B(n,k)$ are known as {\em ballot numbers} and are given by the explicit formula:
\begin{equation}\label{ballot}
B(a,b)=\frac{a-b}{a+b}{a+b\choose a}.
\end{equation}
This expression can be obtained shortly by checking the recurrence \pref{recc}.
We can alternatively use the reflection principle (see \cite{hilton} for a clear presentation), or the cycle lemma ({\it cf.} \cite{cylem}), which will be used in the next sections to obtain a formula in the general case.

The expression \pref{ballot} constitutes a solution to the ballot problem. For this, we use a classical interpretation in terms of paths: we represent a vote for A by an up step, and a vote for B by a down step. The total number of countings of votes, or of paths from $(0,0)$ to $(a+b,a-b)$, is given by the binomial ${a+b\choose a}$. The countings such that A stays ahead of B correspond to paths remaining above the horizontal axis. Their number is given by $B(a,b)$. This implies that the probability asked for at the beginning of this section is $\frac{a-b}{a+b}$.

Of course , we could have used expression \pref{ballot} to prove Proposition \ref{prop1} by a simple computation, but our proof explains more about the combinatorial objects and can be adapted to ternary trees in the next section.

\medskip

It should be mentionned here that our study of multivariate Fuss-Catalan has nothing to do with the many-candidate ballot problem, as considered for example in \cite{zeilb} or \cite{nath}. In particular, the multivariate ballot numbers considered in these papers do not sum to the Fuss-Catalan numbers $C_p(n)=\frac 1 {(p-1)n+1} {pn\choose n}$. Conversely, the numbers studied in the present article do not give any answer to the generalized ballot problem.

\noindent {\bf Remark 1.3.}
Our Catalan array presents similarities with Riordan arrays, but it is not a Riordan array. It may be useful to explicate this point. We recall ({\it cf.} \cite{riordan}) that a Riordan array $M=(m_{i,j})\in\C^{\N\times\N}$ is defined with respect to 2 generating functions 
$$g(x)=\sum g_n x^n \ \ \ {\rm and}\ \ \ f(x)=\sum f_n x^n$$
and is such that
$$M_j(x)=\sum_{n\ge 0} m_{n,j} x^n=g(x){f(x)}^j.$$
We easily observe that a Riordan array with the first two columns $M_0(x)$ and $M_1(x)$ of our Catalan array is relative to $g(x)=\frac{1}{1-x}$ and $f(x)=\frac{x}{1-x}$, which gives the Pascal triangle.

In fact, the Riordan array relative to $g(x)={\bf C}(x)=\sum C(n) x^n$ and $f(x)=x\, {\bf C}(x)$ gives the ballot numbers, but requires knowledge of the Catalan numbers.

\section{Fuss-Catalan tetrahedron and ternary trees}

\subsection{Definitions}

This section, which is the heart of this work, is the study of a 3-dimensional analogue of the Catalan triangle of the previous section. That is we consider exactly the same recurrence, and let the array grow, not in 2, but in 3 dimensions. More precisely, we introduce the sequence $B_3(n,k,l)$ indexed by a positive integer $n$ and nonnegative integers $k$ and $l$, and defined recursively by:
\begin{itemize}
\item $B_3(1,0,0)=1$;
\item $\forall n>1$, $k+l<n$, $B_3(n,k,l)=\sum_{0\le i\le k, 0\le j\le l} B_3(n-1,i,j)$;
\item $\forall k+l\ge n$, $B_3(n,k,l)=0$.
\end{itemize}

Observe that the recursive formula in the second condition is equivalent to: 
\begin{equation}\label{rec}
B_3(n,k,l)=B_3(n-1,k,l)+B_3(n,k-1,l)+B_3(n,k,l-1)-B_3(n,k-1,l-1) 
\end{equation}
and this expression can be used to make some computations lighter, but the presentation above explains more about the generalization of the definition of the ballot numbers $B(n,k)$.

Because of the planar structure of the sheet of paper, we are forced to present the tetrahedron of $B_3(n,k,l)$'s by its sections with a given $n$.

\vskip 0.2cm
$
n=1 \longrightarrow 
\left[
\begin{array}{r}
1
\end{array}
\right]
$

$
n=2 \longrightarrow 
\left[
\begin{array}{rr}
1&1\cr
1&
\end{array}
\right]
$

$
n=3 \longrightarrow 
\left[
\begin{array}{rrr}
1&2&2\cr
2&3&\cr
2&&
\end{array}
\right]
$

$
n=4 \longrightarrow 
\left[
\begin{array}{rrrr}
1&3&5&5\cr
3&8&10&\cr
5&10&&\cr
5&&&
\end{array}
\right]
$

$
n=5 \longrightarrow 
\left[
\begin{array}{rrrrr}
1&4&9&14&14\cr
4&15&30&35&\cr
9&30&45&&\cr
14&35&&&\cr
14&&&&
\end{array}
\right]
$
\vskip 0.2cm

It is clear that $B_3(n,k,0)=B_3(n,0,k)=B(n,k)$. The reader may easily check that when we compute $\sum_{k,l} B_3(n,k,l)$, we obtain: $1,\ 3,\ 12,\ 55,\ 273$. These integers are the first terms of the following sequence ({\it cf.} \cite{njas}):
$$C_3(n)=\frac 1 {2n+1} {3n\choose n}.$$

This fact will be proven in Proposition \ref{prop2}. 

\subsection{Combinatorial interpretation}

Fuss\footnote{Nikolai Fuss (Basel, 1755 -- St Petersburg, 1826) helped Euler prepare over 250 articles for publication over a period on about seven years in which he acted as Euler's assistant, and was from 1800 to 1826 permanent secretary to the St Petersburg Academy.}-Catalan numbers ({\it cf.} \cite{FC,hilton}) are given by the formula
\begin{equation}\label{fuss}
C_p(n)=\frac1{(p-1)n+1}{pn\choose n},
\end{equation}
and $C_3(n)$ appear as order-3 Fuss-Catalan numbers. The integers $C_3(n)$ are known \cite{njas} to count {\em ternary trees}, {\it ie.} trees in which every internal node has exactly 3 sons. 

\vskip 0.2cm
{\centerline{\epsffile{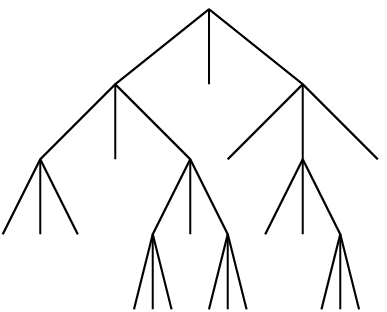}}}

Ternary trees are in bijection with {\em 2-Dyck paths}, which are defined as paths from $(0,0)$ to $(3n,0)$ with steps $(1,1)$ and $(1,-2)$, and remaining above the line $y=0$. The bijection between these objects is the same as in the case of binary trees, {\it ie.} a depth-first search, with the difference that here an internal node is translated into a $(1,-2)$ step. To illustrate this bijection, we give the path corresponding to the previous example of ternary tree:

{\centerline{\epsffile{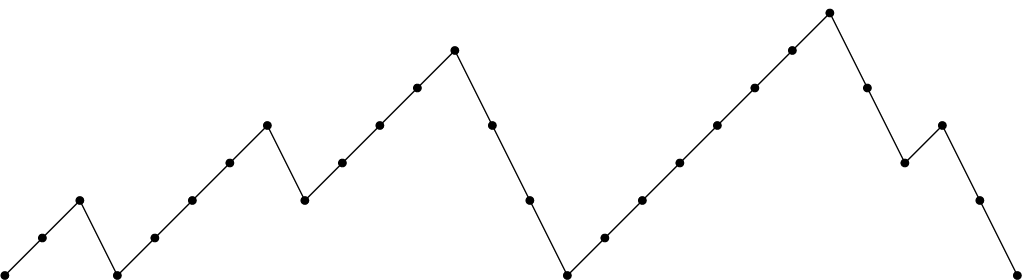}}}

We shall consider these paths with respect to the position of their down steps. The {\em height} of a down step is defined as the height of its end-point. Let $\D_{n,k,l}$ denote the set of 2-Dyck paths of length $3n$, with $k$ down steps at even height and $l$ down steps at odd height, excluding the last seqence of down steps. By definition, the last sequence of down steps is of length $n-k-l$.

\begin{prop}\label{prop2}
We have
$$\forall n>0,\ \ \ \sum_{k,l} B_3(n,k,l)=C_3(n)=\frac 1 {2n+1} {3n\choose n}.$$
Moreover, $B_3(n,k,l)$ is the cardinality of $\D_{n,k,l}$.
\end{prop}

\begin{proof}
Let $k$ and $l$ be fixed. Let us consider an element of $\D_{n,k,l}$. If we cut this path after its $(2n-2)$-th up step, and complete with down steps, we obtain a 2-Dyck path of length $3(n-1)$ (see figure below). It is clear that this path is an element of $\D_{n,i,j}$ for some $i\le k$ and $j\le l$. We can furthermore reconstruct the original path from the truncated one, if we know $k$ and $l$. We only have to delete the last sequence of down steps (here the dashed line), to draw $k-i$ down steps, one up step, $l-j$ down steps, one up step, and to complete with down steps. This gives a bijection from $\D_{n,k,l}$ to $\cup_{0\le i\le k,0\le j\le l}\D_{n-1,i,j}$, which implies Proposition \ref{prop2}.

{\centerline{\epsffile{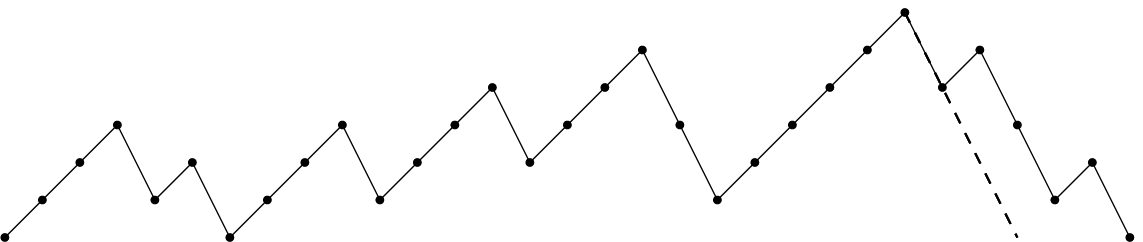}}}

\end{proof}

\noindent {\bf Remark 2.2.}
It is interesting to translate the bi-statistic introduced on 2-Dyck paths to the case of ternary trees. As previously, we consider the depth-first search of the tree, and shall not consider the last right string. We define $\T_{n,k,l}$ as the set of ternary trees with $n$ internal nodes, $k$ of them being encountered in the search after an even number of leaves and $l$ after and odd number of leaves. By the bijection between trees and paths, and Proposition \ref{prop2}, we have that the cardinality of $T_{n,k,l}$ is $B_3(n,k,l)$.

\noindent {\bf Remark 2.3.}
It is clear from the definition that:
$$B_3(n,k,l)=B_3(n,l,k).$$
But this fact is not obvious when considering trees or paths, since the statistics defined are not clearly symmetric. To explain this, we can introduce an involution on the set of ternary trees which sends an element of $\T_{n,k,l}$ to $\T_{n,l,k}$. To do this, we can exchange for each node of the last right string its left and its middle son, as in the following picture. Since the number of leaves of a ternary tree is odd, every ``even'' node becomes an odd one, and conversely.

{\centerline{\epsffile{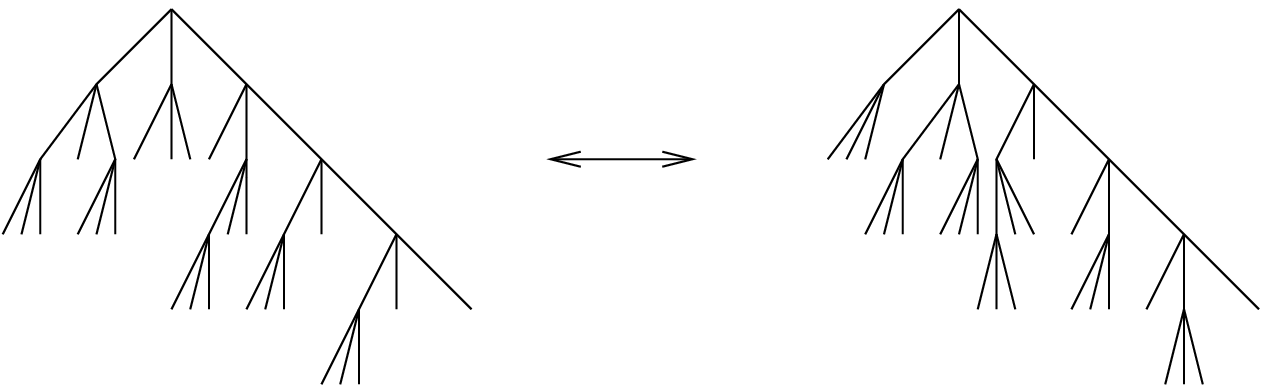}}}

\subsection{Explicit formula}

Now a natural question is to obtain explicit formulas for the $B_3(n,k,l)$. The answer is given by the following proposition.

\begin{prop}\label{prop2.1}
The integers $B_3(n,k,l)$ are given by
\begin{equation}\label{prop2.1eq}
B_3(n,k,l)={n+k\choose k}{n+l-1\choose l}\frac{n-k-l}{n+k}
\end{equation}
\end{prop}

\begin{proof}
We use a combinatorial method to enumerate $\D_{n,k,l}$. 
The method is a variation of the cycle lemma \cite{cylem} (called ``penetrating analysis'' in \cite{gopal}).

If we forget the condition of ``positivity'' ({\it ie.} the path remains above the line $y=0$), and cut the last down sequence, a path consists in:
\begin{itemize}
\item at even height: $n$ up steps, and $k$ down steps;
\item at odd height: $n$ up steps and $l$ down steps.
\end{itemize}
An important remark is to remember that an element of $\D_{n,k,l}$ has an up step just before the last sequence of down steps! 
If we suppose that all these steps are distinguished, we obtain:
\begin{itemize}
\item $\frac{(n+k)!}{n!k!}$ choices for even places;
\item $\frac{(n-1+l)!}{(n-1)!l!}$ choices for odd places (we cannot put any odd down step after the last odd step).
\end{itemize}
Now we group the paths which are ``even permutations'' of a given path $P$. By even permutations, we mean cycle permutations which preserve the parity of the height of the steps. 

We want to prove that the proportion of elements of $\D_{n,k,l}$ in any even orbit (\ie in any orbit under even permutations) is given by $\frac{n-k-l}{n+k}$.

We suppose first that the path $P$ is acyclic (as a word): $P$ cannot be written as $P=U^p$ with $U$ a word in the two letters $(1,1)$ and $(1,-2)$. It is clear that such a path gives $n+k$ different even permutations.
Now we have to keep only those which give elements of $\D_{n,k,l}$. 
To do this we consider the concatenation of $P$ and $P'$, which is a duplicate of $P$.

\centerline{\epsffile{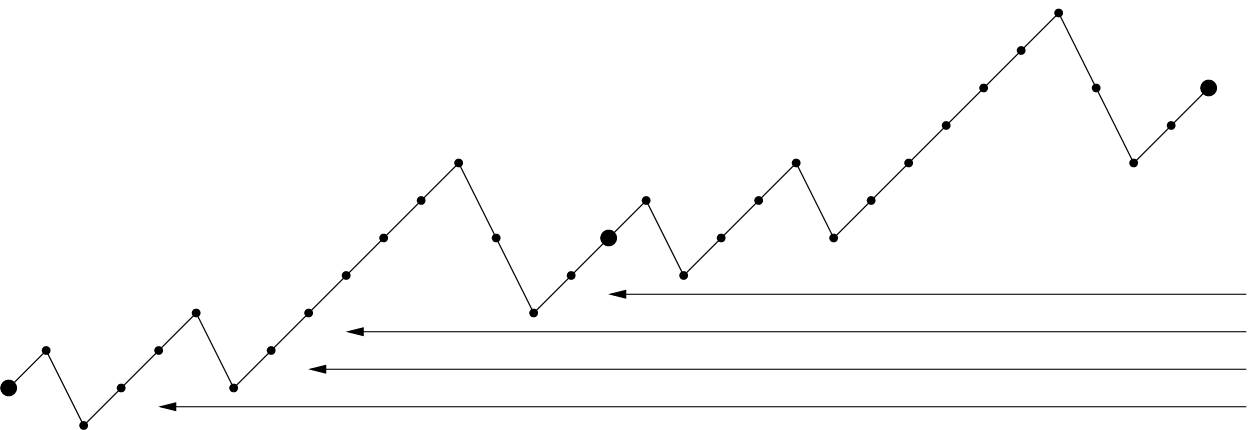}}

The cyclic permutations of $P$ (not necessarily even) are the subpaths of $P+P'$ of horizontal length $2n+k+l$. The number of such paths that remain above the horizontal axis, and end with an up step, is the number of (up) steps of $P$ in the light of an horizontal light source coming from the right. The only transformation is to put the illuminated up step at the end of the path. The number of illuminated up steps is $2n-2k-2l$, since every down step puts two up steps in the shadow. Among these $2n-2k-2l$ permutations, only half are even (observe that the set of heights of the illuminated steps is an interval). Thus $n-k-l$ paths among the $n+k$ elements of the orbit are in $\D_{n,k,l}$.

Now we observe that if $P$ is $p$-cyclic, then its orbit has $p$ times less elements, and we obtain $p$ times fewer different paths, whence the proportion of elements of elements of $\D_{n,k,l}$ in this orbit is:
$$\frac{(n-k-l)/p}{(n+k)/p}=\frac{n-k-l}{n+k}.$$

Finally, we obtain that the cardinality of $\D_{n,k,l}$ is 
$$\frac{(n+k)!}{n!k!}\frac{(n-1+l)!}{(n-1)!l!}\frac{n-k-l}{n+k}$$
which was to be proved. 

\end{proof}

\noindent {\bf Remark 2.5.}
The equation \pref{prop2.1eq} is of course symmetric in $k$ and $l$:
$${n+k\choose k}{n+l-1\choose l}\frac{n-k-l}{n+k}={n+k-1\choose k}{n+l-1\choose l}\frac{n-k-l}{n}.$$

\noindent {\bf Remark 2.6.}
I made the choice to present a combinatorial proof of Proposition \ref{prop2.1}. The interest is to show how these formulas are obtained, and to allow easy generalizations ({\it cf.} the next section). It is also possible to check directly the recurrence \pref{rec}.

\subsection{Generating function}

Let $F$ denote the following generating series:
\begin{equation}\label{F}
F(t,x,y)=1+\sum_{0\le k+l<n} B_3(n,k,l) t^n x^k y^l.
\end{equation}
In this expression, the term ``1'' corresponds to the {\em empty} 2-Dyck path. Thus $F$ is the generating series of 2-Dyck paths with respect to their length (variable $t$), their number of down steps at even height --excluding the last down sequence-- (variable $x$), and their number of down steps at odd height (variable $y$).

We also introduce 
\begin{equation}\label{G}
G(t,x,y)=1+\sum_{0\le k+l<n} B_3(n,k,l) t^n x^{n-l} y^l,
\end{equation}
{\it ie.} $G$ is the generating function of 2-Dyck paths with respect to their length, to the number of down steps at even height -including the last down sequence- and to the number of down steps at odd height.

To obtain equations for $F$ and $G$, we decompose a non-empty 2-Dyck path by looking at two points: $\alpha$, defined as the last return to the axis (except the final point $(3n,0)$), and $\beta$ defined as the last point at height 1 after $\alpha$. This gives the following decomposition of a 2-Dyck path
$$P=P_1\ (up)\ P_2\ (up)\ P_3\ (down),$$
with $P_1$, $P_2$, $P_3$ any 2-Dyck path (maybe empty).

\vskip 0.2cm
\centerline{\epsffile{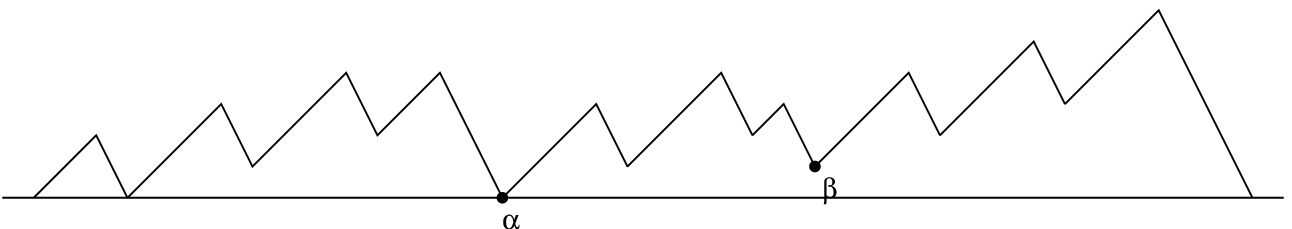}}
\vskip0.2cm

By observing that the up steps after $\alpha$ and $\beta$ change the parity of the height,
this gives the two following equations :

\begin{equation}
F(t,x,y)=1+G(t,x,y)\times G(t,y,x)\times t.F(t,x,y)
\end{equation}
and in the same way
\begin{equation}\label{Geq}
G(t,x,y)=1+G(t,x,y)\times G(t,y,x)\times t.G(t,x,y).x.
\end{equation}

By permuting the variables $x$ and $y$ in \pref{Geq}, we obtain
\begin{equation}\label{Geq2}
G(t,y,x)=1+G(t,y,x)\times G(t,x,y)\times t.G(t,y,x).y
\end{equation}
and we can eliminate $G(t,y,x)$ from \pref{Geq} and \pref{Geq2} to obtain the following result.

\begin{prop}
The generating function $F$ of the $B_3$'s is given by:
$$F(t,x,y)=\frac{1}{1-tG(t,x,y)G(t,y,x)}$$
where $G(t,x,y)$ is a solution of the algebraic equation
$$tx^2G^3+(y-x)G^2+(x-2y)G+y=0.$$
\end{prop}

\vskip 0.3cm
An alternate approach to the generating function is to use formula \ref{prop2.1eq} to obtain what MacMahon called a ``redundant generating function `` ({\it cf.} \cite{gessel}), since it contains terms other than those which are combinatorially significant.

To do this we extend the recursive definition of $B_3(n,k,l)$ as follows: we define $B'_3(n,k,l)$ for integers $n>0$, and  $k,l\ge0$ by:
\begin{equation}\label{prime}
B'_3(n,k,l)={n+k-1\choose k}{n+l-1\choose l}\frac{n-k-l}{n}.
\end{equation}

Of course when $k+l>n$, the integer $B'_3(n,k,l)$ is negative.

As an example, here is the ``section'' of the array of $B'_3(n,k,l)$ with $n=3$.

\vskip 0.2cm
$\left[\begin{array}{cccccc}
1&2&2&0&-5&\cdots\cr
2&3&0&-10&-30&\cdots\cr
2&0&-12&-40&-90&\cdots\cr
0&-10&-40&-100&-200&\cdots\cr
-5&-30&-90&-200&-375&\cdots\cr
\vdots&\vdots&\vdots&\vdots&\vdots&\ddots 
\end{array}
\right]
$
\vskip 0.2cm

The equation \pref{prime} is equivalent to the following recursive definition: 
 $\forall n\le 0$ or $k,l<0,\ B'_3(n,k,l)=0$ and 
$$B'_3(n,k,l)=B'_3(n-1,k,l)+B'_3(n,k-1,l)+B'_3(n,k,l-1)-B'_ 3(n,k-1,l-1)+c(n,k,l)$$
$${\rm where}\ \ \ \ c(n,k,l)=\left\{
\begin{array}{rl}
+1&{\rm if }(n,k,l)=(1,0,0)\cr
-1&{\rm if }(n,k,l)=(0,1,0)\cr
-1&{\rm if }(n,k,l)=(0,0,1)\cr
+2&{\rm if }(n,k,l)=(0,1,1)\cr
0&{\rm otherwise.}
\end{array}
\right.
$$

The interest of this presentation is to give a simple (rational!) generating series. Indeed, it is quite simple to deduce from the recursive definition of $B'_3(n,k,l)$ that:
\begin{equation}\label{F'}
\sum_{n,k,l}B'_3(n,k,l)t^nx^ky^l=\frac{t-x-y+2xy}{1-t-x-y+xy}.
\end{equation}

This formula then encodes the $B_3(n,k,l)$ since $\forall k+l<n,\ B_3(n,k,l)=B'_3(n,k,l)$.

\section{Fuss-Catalan $p$-simplex and $p$-ary trees}

The aim of this final section is to present an extension of the results of Section 2 to $p$-dimensional recursive sequences. In the same spirit, we define the sequence $B_p(n,k_1,k_2,\dots,k_{p-1})$ by the recurrence:
\begin{itemize}
\item $B_p(1,0,0,\dots,0)=1$;
\item $\forall n> 1$ and $0\le k_1+k_2+\cdots+k_{p-1}< n$, 
$$B_p(n,k_1,k_2,\dots,k_{p-1})=\sum_{0\le i_1\le k_1,\dots,0\le i_{p-1}\le k_{p-1}} B(n-1,i_1,i_2,\dots,i_{p-1});$$
\item $\forall k_1+k_2+\cdots+k_{p-1}\ge n$, $B_p(n,k_1,k_2,\dots,k_{p-1})=0$.
\end{itemize}

\noindent
Every result of Section 2 extends to general $p$. We shall only give the main results, since the proofs are straightforward generalizations of the proofs in the previous section.

\begin{prop}\label{prop3}
$$\sum_{k_1,\dots,k_{p-1}} B_p(n,k_1,\dots,k_{p-1})=C_p(n)=\frac 1 {(p-1)n+1} {pn\choose n}$$
\end{prop}

The integers $C_p(n)$ are order-$p$ Fuss-Catalan numbers and enumerate $p$-ary trees, or alternatively $p$-Dyck paths (the down steps are $(1,-p)$). In this general case, the recursive definition of $B_p(n,k_1,k_2,\dots,k_{p-1})$ gives rise to $p-1$ statistics on trees and paths analogous to those defined in section 3.

\noindent {\bf Remark 3.2.}
By the same method as in the previous section, it is possible to obtain an explicit formula for these multivariate Fuss-Catalan numbers:
$$B_p(n,k_1,k_2,\dots,k_{p-1})=\left(\prod_{i=1}^{p-1}{n+k_i-1\choose k_i}\right)\frac{n-\sum_{i=1}^{p-1}k_i}{n}.$$


\vskip 0.3cm

\noindent
{\bf\large Comment.} The numbers $B_3(n,k,l)$ first arose in a question of algebraic combinatorics ({\it cf.} \cite{aval}). Let $\I_n$ be the ideal generated by $B$-quasisymmetric polynomials in the $2n$ variables $x_1,\dots,x_n$ and $y_1,\dots,y_n$ ({\it cf.} \cite{BH}) without constant term. We denote by $\QQ_n$ the quotient $\Q[x_1,\dots,x_n,y_1,\dots,y_n]/\I_n$ and by $\QQ_n^{k,l}$ the bihomogeneous component of $\QQ_n$ of degree $k$ in $x_1,\dots,x_n$ and degree $l$ in $y_1,\dots,y_n$. It is proven in \cite{aval} that:
$$\dim \QQ_n^{k,l} = B_3(n,k,l) = {n+k-1\choose k}{n+l-1\choose l}\frac{n-k-l}{n}.$$


\vskip 0.3cm

\noindent
{\bf\large Acknowledgement.} The author is very grateful to referees (and in particular to the anonymous ``Referee 3'') for valuable remarks and corrections.



\vskip 0.3cm

\begin{thebibliography}{10}

\bibitem{andre} 
 {\sc D. Andr\'e}, 
 {\em Solution directe du probl\`eme r\'esolu par M. Bertrand},
 Comptes Rendus Acad. Sci. Paris, {\bf 105} (1887), 436--437.

\bibitem{aval} 
 {\sc J.-C. Aval}, 
 {\em Ideals and quotients of $B$-quasisymmetric functions},
 Sem. Lothar. Comb. 54 (2006), B54d (13 pages).

\bibitem{eco} 
 {\sc E. Barcucci, A. Del Lungo, E. Pergola, R. Pinzani}, 
 {\em ECO: a methodology for the enumeration combinatorial objects}, 
 J. Difference Equations Appl. 5 (1999) 435--490.

\bibitem{BH}
{\sc P. Baumann and C. Hohlweg}, 
{\em A Solomon theory for the wreath product $G\wr{\mathfrak S}_n$}, preprint.

\bibitem{bertrand} 
 {\sc J. Bertrand}, 
 {\em Calcul des probabilit\'es. Solution d'un probl\`eme},
 Comptes Rendus Acad. Sci. Paris, {\bf 105} (1887), 369.

\bibitem{cylem} 
 {\sc N. Dershowitz and S. Zaks}, 
 {\em The Cycle Lemma and Some Applications},
 Europ. J. Combinatorics, {\bf 11} (1990), 35--40.

\bibitem{euler1}
  {\sc L. Euler}, 
  {\em Leonhard Euler und Christian Goldbach, Briefwechsel 1729-1764}, 
    Juskevic, A. P., Winter, E. (eds.), Akademie Verlag, Berlin, 1965.

\bibitem{euler2} 
  \sc {L. Euler}, 
  {\em Novi Commentarii Academiae Scientarium Imperialis Petropolitanque 7},
  (1758-59), pp. 9-28.

\bibitem{FC}
 {\sc S. Fomin and N. Reading},
 {\em Generalized cluster complexes and Coxeter combinatorics},
 In preparation.

\bibitem{gessel} 
 {\sc I. M. Gessel}, 
 {\em Super Ballot Numbers},
 J. Symbolic Computation, {\bf 14} (1992), 179--194.

\bibitem{hilton} 
 {\sc P. Hilton and J. Pedersen}, 
 {\em Catalan Numbers, their Generalizations, and their Uses},
 Math. Intelligencer, {\bf 13} (1001), 64--75.

\bibitem{gopal} 
 {\sc S. G. Mohanty}, 
 {\em Lattice paths counting and application},
 New York and London: Academic Press (1979).

\bibitem{nath} 
 {\sc M. B. Nathanson}, 
 {\em Ballot numbers, alternating products and the Erd\"os-Heilbronn conjecture},
 in: R. L. Graham, J. Nesetril (Eds.), The Mathematics of Paul Erd\"os, Springer-Verlag, Berlin, 1997, pp 199-217.

\bibitem{riordan} 
 {\sc L. W. Shapiro, S. Getu, W.-J. Woan and L. C. Woodson}, 
 {\em The Riordan Group},
 Discrete Appl. Math., {\bf 34} (1991), 229--239.

\bibitem{njas} 
 {\sc N. J. A. Sloane}, 
 {\em A Handbook of integer sequences},
 New York and London: Academic Press (1973).

\bibitem{stanley} 
 {\sc R. Stanley}, 
 {\em Enumerative Combinatorics, volume 2},
 no.~62 in
 Cambridge 
    Studies in Advanced Mathematics, Cambridge University Press, 1999.

\bibitem{stanweb} 
 {\sc R. Stanley}, 
 {\em Catalan Addendum},
 (web) {\tt http://www-math.mit.edu/~rstan/ec/catadd.pdf}.

\bibitem{zeilb} 
 {\sc D. Zeilberger}, 
 {\em Andr\'e's reflection proof generalized to the many-candidate ballot problem},
 Disc. Math., {\bf 44} (1983), 325-326.

\end{thebibliography}
\end{document}